\documentclass[12pt, reqno]{amsart}
\usepackage{color}
\usepackage{amsmath}
\usepackage{amssymb}
\usepackage{amsfonts}
\usepackage{mathrsfs}
\usepackage{amsbsy}
\usepackage{relsize}
\usepackage[colorlinks, linkcolor=red, citecolor=blue, urlcolor=blue, pagebackref, hypertexnames=false]{hyperref}
\usepackage[hyperpageref]{backref}

\usepackage[lmargin=3cm,rmargin=3cm,tmargin=3cm,bmargin=3cm]{geometry}
\newcommand{\beq}{\begin{eqnarray}}
\newcommand{\eeq}{\end{eqnarray}}
\newcommand{\bq}{\begin{equation}}
\newcommand{\eq}{\end{equation}}
\newcommand{\beqn}{\begin{eqnarray*}}
\newcommand{\eeqn}{\end{eqnarray*}}

\newcommand{\R}{\ensuremath{\mathbb{R}}}

\newcommand{\vertiii}[1]{{\vert\kern-0.25ex\vert\kern-0.25ex\vert #1
    \vert\kern-0.25ex\vert\kern-0.25ex\vert}}

\newcommand{\ignore}[1]{}

\newtheorem{definition}{Definition}[section]
\newtheorem{proposition}{Proposition}[section]
\newtheorem{theorem}{Theorem}[section]
\newtheorem{remark}{Remark}[section]

\title[Hardy-H\'enon parabolic equations with fractional noise]{Well-posedness for Hardy-H\'enon parabolic equations with fractional Brownian noise}

\author[M. Majdoub, E. Mliki]{Mohamed Majdoub \& Ezzedine Mliki}
\address{Department of Mathematics, College of Science, Imam Abdulrahman Bin Faisal University, P. O. Box 1982, Dammam, Saudi Arabia}
\address{Basic and Applied Scientific Research Center, Imam Abdulrahman Bin Faisal University, P.O. Box 1982, 31441, Dammam, Saudi Arabia}
\email{\sl mmajdoub@iau.edu.sa}
\email{\sl ermliki@iau.edu.sa}

\begin{document}
\begin{abstract} We study the Hardy-H\'enon parabolic equations on $ \mathbb{R}^{N}$ ($N$=2 or 3) under the effect of an additive fractional Brownian noise with Hurst parameter $H>\max\left(1/2, N/4\right).$ We show local existence and uniqueness of a mid $L^{q}$-solution under suitable assumptions on $q$.
\end{abstract}




\date{\today}

\maketitle

\section{Introduction}
We consider the following Cauchy problem
\begin{equation}
\label{main}
\left\{
\begin{matrix}
\partial_tu(t)-\Delta u(t)=|x|^{-\gamma} |u(t)|^{p-1}u(t)+\partial_t B^H(t),\quad t>0\\
u(0)= u_0,\\
\end{matrix}
\right.
\end{equation}
where $p>1$, $x\in\R^N$, $u_0\in L^q(\R^N)$, $\gamma\geq 0$, and the random forcing $B^H$ is the fractional Brownian motion defined on some complete probability space $(\Omega, \mathcal{F}, \mathbb{P})$ with Hurst parameter $H\in (0,1)$. When $H=1/2$, $B^{1/2}$ is the Wiener process.\\

 As is a standard practice, we study \eqref{main} via the associated integral equation:
\begin{equation}
\label{integral}
u(t)= {\rm e}^{t\Delta}u_{0}+\int_{0}^{t}{\rm e}^{(t-s)\Delta}\,\left(|x|^{-\gamma}|u(s)|^{p-1}u(s)\right)\,ds+\int_0^t{\rm e}^{(t-s)\Delta}\,dB^H(s),
\end{equation}
where ${\rm e}^{t\Delta}$ is the linear heat semi-group.

The initial value problem \eqref{main} without fractional noise has attracted considerable attention in the mathematical community and the well-posedness theory in the Lebesgue spaces. The first works in this direction are due to Weissler \cite{W1, W2, W3} when $\gamma=0$.
Removing the fractional noise and assuming $\gamma=0$, we see that the equation enjoys an interesting property of scaling invariance
\begin{equation}
\label{scaling}
u_\lambda(t,x):=\lambda^{{2\over p-1}}\,u(\lambda^2 t, \lambda x), \quad \lambda>0.
\end{equation}
Note that the $L^q(\R^N)$ norm is invariant under this scaling if and only if
$$
q=q_c=\frac{N(p-1)}{2}.
$$
This leads the following classification with regards to existence and uniqueness:
\begin{itemize}
\item [{\sc Case 1.}] If $q\geq q_c$ and $q>1$ or $q>q_c$ and $q\geq 1$, Weissler in \cite{W2} proved the existence of a unique local solution $u\in {\mathcal C}([0,T); L^q(\R^N))\cap L_{loc}^\infty(]0,T]; L^\infty(\R^N))$. Later on, Brezis-Cazenave \cite{Br} proved the unconditional uniqueness of Weissler's solutions.\footnote{Uniqueness in the natural space where solutions exist, namely ${ C}([0,T); L^q)$.}
\item [{\sc Case 2.}] If $q<q_c$, there are indications that there exists no (local) solution in any reasonable weak sense. See \cite{Br, W2, W3}. Moreover, it is known that uniqueness is lost for the initial data $u_0=0$ and for $1+\frac{1}{N}<p<\frac{N+2}{N-2}$. See \cite{HW}.
\end{itemize}

Concerning blow-up there are two effective techniques which has been employed to prove non-existence of global solutions: the concavity method (\cite{Lev}) and the eigenfunction method (\cite{Kap}). The later one was firstly used for bounded domains but can be adapted to the whole space $\R^N$. The concavity method and its variants were used in the study of many nonlinear evolution partial differential equations (see e.g. \cite{GP1, GP2, PS}).

Recently, the case of exponential nonlinearity was considered and results on global existence, blow-up and decay estimates was obtained (see \cite{IJMS, Ioku, MT, MT1, IRT, RT}). See also \cite{MOT} for the biharmonic heat equation. Let us mention that the well posedness in Sobolev and Besov spaces was investigated in \cite{Rib, MI}.\\

Note that the case $\gamma>0$ was investigated in \cite{BTW}. We refer the reader to \cite{QS} and references therein for more properties and information on nonlinear heat equations.\\

Nowadays, the study of stochastic partial differential equation driven by fractional noise has attracted considerable attention in the mathematical community, motivated by theoretical reasons and also by its applications in physics, biology, hydrology, and other sciences.  A special interest has been attached to the well-posedness for semilinear stochastic parabolic equations driven by an infinite-dimensional fractional noise. See for instance \cite{DDM, MS, ND, SV}. Other type of noises have been also considered in \cite{BT, BNS,CM}.\\

Here we consider the semilinear Hardy-H\'enon equation with fractional Brownian noise \eqref{main}. Our main interest is to solve \eqref{main} locally in time for initial data in Lebesgue space $L^q$. We first define what we call a mild solution of \eqref{main}.
\begin{definition}
\label{mild}
A measurable function $u : \Omega\times [0,T]\to L^q$ is a mild solution of the integral equation \eqref{integral} if
\begin{itemize}
\item[i)] $u$ satisfies \eqref{integral} with probability one.
\item[ii)] $u\in C([0,T]; L^q)$.
\end{itemize}
\end{definition}

We have obtained the following local well-posedness result.

\begin{theorem}
\label{LWP}
Let $N=2, 3$. Assume that
\begin{equation}
\label{gamma}
0< \gamma<2,
\end{equation}
\begin{equation}
\label{H}
\max\left(1/2, N/4\right)<H<1,
\end{equation}
and
\begin{equation}
\label{Q}
\max\left(\frac{Np}{N-\gamma}, \frac{N(p-1)}{2-\gamma}\right)<q<\infty.
\end{equation}

Suppose that $u_0\in L^q(\R^N)$. Then there exists $T=T(\|u_0\|_q)>0$ such that problem \eqref{main} possesses a unique mild solution in $[0,T]$.
\end{theorem}
\begin{remark}
The Cauchy problem was studied in \cite{CO} for $\gamma=0$. The local existence was proved under restrictive assumptions.
\end{remark}

The paper is organized as follows. In Section 2, we recall some preliminaries needed in the paper such as fractional Brownian motion and smoothing effect for the heat semi-group. The third section is devoted to the proof of our main result Theorem \ref{LWP}.

\section{Preliminaries}

\subsection{Fractional Brownian motion}
In the following $(\Omega, \mathcal{F}, \mathbb{P})$ stands for a complete probability space.
Fractional Brownian motion was originally defined and study by Kolmogorov \cite{Kol} within a Hilbert space framework. Consider $[0, T]$ a time interval with arbitrary fixed horizon $T,$ a  fractional Braownian motion with Hurst parameter $H\in (0, 1)$ is a centered Gaussian process $B^{H}$ with covariance
\begin{equation}\label{cov}
R(s, t):=E(B^{H}(t)B^{H}(s))= \frac{1}{2} \left(t^{2H}+s^{2H}- |t-s|^{2H}\right)
\end{equation}
where $s,t\in [0, T].$\\
If $H=\frac{1}{2}$ then $B^{\frac{1}{2}}$ is standard Brownian motion.
The fractional Brownian motion ({\tt fBm}) can also be defined as the only self-similar Gaussian process with stationary increments
\begin{equation*}
E[(B^{H}(t)-B^{H}(s))^{2}]=|t-s|^{2H}
\end{equation*}
and $H$-self similar
\begin{equation*}
(\frac{1}{c^{H}}B^{H}(ct), \, t\geq0)=^{d}(B^{H}(t), \, t\geq0)
\end{equation*}
for all $c>0.$
Moreover the process $B^{H}$ has the following Wienner integral representation:
\begin{equation}\label{BH}
B^{H}(t)=\int_{0}^{t}K^{H}(t, s)\, dW(s)
\end{equation}
where $W=\{W(t);\, t\in [0, T]\}$ is a Wiener process, and $K^{H}(t, s)$ is the kernel given by
 \begin{equation}\label{KH}
 K^{H}(t, s)=c_{H}(t-s)^{H-\frac{1}{2}}+c_{H} (\frac{1}{2}-H)\int_{s}^{t}(r-s)^{H-\frac{3}{2}}(1-(1+(\frac{s}{r}))^{\frac{1}{2}-H})\, dr
 \end{equation}
 $c_{H}$ is a constant given by
 \begin{equation*}
 c_{H}=\left(\frac{2H\Gamma(\frac{3}{2}-H)}{\Gamma(H+\frac{1}{2})\Gamma(2-2H)}\right)^{\frac{1}{2}}
   \end{equation*}
 From \eqref{KH} we obtain
  \begin{equation*}
  \frac{\partial K^{H}}{\partial t}(t, s)=c_{H}(H-\frac{1}{2})\left(\frac{s}{t}\right))^{\frac{1}{2}-H}(t-s)^{H-\frac{3}{2}}
   \end{equation*}
   Notice that if $H>\frac{1}{2}$ then the kernel $K_{H}(t,s)$ is regular and has the simpler expression
   \begin{equation*}
   K^{H}(t,s)=c_{H} s^{\frac{1}{2}-H}\int_{s}^{t}(r-s)^{H-\frac{3}{2}}r^{H-\frac{1}{2}}\,dr
   \end{equation*}
We will denote by $\varepsilon_{H}$ the linear space of step functions on $[0, T]$ of the form
 \begin{equation}\label{pphi}
 \varphi(t)=\sum_{i=1}^{n}a_{i}1_{(t_{i}, t_{i+1}]}(t)
  \end{equation}
  where $t_{1}, ... , t_{n}\in [0,T],$ $n\in \mathbb{N},$ $a_{i}\in \mathbb{R}$ and by $\mathcal{H}$ the closure of $\varepsilon_{H}$ with respect to the scalar product
   \begin{equation*}
\langle1_{[0, t]}, \, 1_{[0, s]}\rangle_{\mathcal{H}}=R(t,s)
  \end{equation*}
  For $\varphi \in\varepsilon_{H}$ of the form \eqref{pphi}  we define its Wiener integral with respect to the fractional Brownian motion as
  \begin{equation*}
  \int_{0}^{T}\varphi_{s}\,dB^{H}(s)=\sum_{ i=1}^{n}a_{i}(B^{H}_{t_{i+1}}-B^{H}_{t_{i}})
  \end{equation*}
  Obviously, the mapping
    \begin{equation*}
     \varphi=\sum_{i=1}^{n}a_{i}1_{(t_{i}, t_{i+1}]}\longrightarrow  \int_{0}^{T}\varphi_{s}\,dB^{H}(s)
      \end{equation*}
      is an isometry between $\varepsilon_{H}$ and the linear space $span\{B_{t}^{H},\, t\in \mathbb{R}\}$ viewed as a subspace of $L^{2}(\Omega).$ The image on an element $\Phi \in \mathcal{H}$ by this isometry is called the Wiener integral of $\Phi$ with respect to $B^{H}$.
      For every $s<T,$ let us consider the operator $K^{\ast}$ in $L^{2}[0, T]$
          \begin{equation*}
          (K^{\ast}_{T}\varphi)(s)=K(T,s)\varphi(s)+\int_{s}^{T}(\varphi(r)-\varphi(s))\frac{\partial K}{\partial r}(r,s)\, dr
          \end{equation*}
      When $H>\frac{1}{2}$, the operator $K^{\ast}$ has the simpler expression
      \begin{equation*}
                (K^{\ast}_{T}\varphi)(s)=\int_{s}^{T}\varphi(r)\frac{\partial K}{\partial r}(r,s)\, dr
          \end{equation*}
    We refer to \cite{AMN} for the proof of the fact that $K^{\ast}$ is an isometry between $\mathcal{H}$ and $L^{2}[0, T]$ and, as a consequence, we will have the following relationship between the Wiener process $W$
     \begin{equation*}
     \int_{0}^{t}\varphi(s)\,dB^{H}(s)=\int_{0}^{t}(K^{\ast}_{t}\varphi)(s)\, dW(s)
     \end{equation*}
     for every $t\in [0, T]$ and $\varphi 1_{[0, t]}\in \mathcal{H}$ if and only if $K^{\ast}\varphi \in L^{2}[0, T].$ We also recall that, if $\phi, \chi \in \mathcal{H}$ are such that $\int_{0}^{T}\int_{0}^{T}|\phi(s)||\chi(t)||t-s|^{2H-2}\, ds dt<\infty,$ their scalar product in $\mathcal{H}$ is given by \begin{equation}\label{intW}
     \langle\phi, \chi\rangle_{H}= H(2H-1)\int_{0}^{T}\int_{0}^{T}\phi(s)\chi(t)|t-s|^{2H-2}\, ds dt
       \end{equation}
       Note that  in the general, the existence of the right-hand side of  \eqref{intW} requires careful justification (see \cite{ND}). As we will work only with Wiener integral over Hilbert space, we point out that if $X$ is a Hilbert space and $u\in L^{2}([0, T]; X)$, is a deterministic function, the relation \eqref{intW} holds, and the right-hand side being well defined in $L^{2}(\Omega, X)$ if $K^{\ast}u$ belongs $L^{2}([0, T]\times X).$
          \subsection{Cylindrical fractional Brownian motion}
          As in \cite{DDM}, we define the standard cylindrical fractional Brownian motion in $X$ as the formal series
          \begin{equation}\label{seri}
          B^{H}(t)=\sum_{n=0}^{\infty}e_{n}b_{n}^{H}(t)
          \end{equation}
          where $\{e_{n},\, n\in \mathbb{N}\}$ is a complete orthonormal basis in $X$, and  $b_{n}^{H}$ is a one dimensional {\tt fBm}. It is well known that the infinite series \eqref{seri} does not converge in $L^{2}(\mathbb{P}),$ hence $B^{H}(t)$ is not a well-defined $X$-valued random variable. Nevertheless, for every Hilbert space ${\mathcal N}$ such that $X\hookrightarrow {\mathcal N}$, the linear embedding is Hilbert-Schmidt operator, therefore, the series \eqref{seri} defines a ${\mathcal N}$-valued random variable and $\{B^{H}(t), \, t\geq0\}$ is a ${\mathcal N}$-valued $Id-fBm$.\\
          Following the approach for cylindrical Brownian motion introduced in \cite{DZ}, it is possible to define a stochastic integral of the form
            \begin{equation}\label{stint}
         \int_{0} ^{T}f(t)B^{H}(t)
          \end{equation}
          where $f:[0, T]\mapsto \mathcal{L}(X, Y)$ and $Y$ is another real and separable Hilbert space, and the \eqref{stint} is a $Y$-valued random variable that is independent of choice of ${\mathcal N}.$\\

          Let $f$ be a deterministic function with values in $\mathcal{L}_{2}(X, Y)$, the space of Hilbert-Schmidt operators from $X$ to $Y.$ We consider the following assumptions on $f.$

          \begin{itemize}
\item[i)] For each $x\in X$, $f(.)x\in L^{p}([0, T]; Y),$ for $p>\frac{1}{H}.$
\item[ii)] $\alpha_{H}\int_{0}^{T}\int_{0}^{T}|f(s)|_{\mathcal{L}_{2}(X, Y)}|f(t)|_{\mathcal{L}_{2}(X, Y)}|s-t|^{2H-2}\,ds dt<\infty.$
\end{itemize}
The stochastic integral \eqref{stint} is defined as

              \begin{equation}\label{sumint}
         \int_{0} ^{T}f(t)B^{H}(t):=\sum_{n=1}^{\infty}\int_{0}^{t} f(s) e_{n}\, db_{n}^{H}(s)=\sum_{n=1}^{\infty}\int_{0}^{t}(K_{H}^{\ast}fe_{n})\, db_{n}(s),
        \end{equation}
          where $b_{n}$ is the standard Brownian motion linked to ({\tt fBm}) $b_{n}^{H}$ via the representation formula \eqref{BH}. Since $fe_{n}\in L^{2}([0, T]; Y)$ for each $n\in \mathbb{N}$, the variables $\{\int_{0}^{t}fe_{n}\, db_{n}^{H}\}$ are mutually independent (see \cite{DDM}).
          The series \eqref{sumint} is finite

          \begin{equation}\label{fiseri}
      \sum_{n}\|K^{\ast}_{H} (fe_{n})\|^{2}\, db_{n}(s)=\sum_{n}\parallel\parallel fe_{n}\parallel_{\mathcal{H}}\parallel_{X}^{2}<\infty
          \end{equation}
          If we consider $X=Y=\mathcal{H},$ we have

             \begin{eqnarray*}
         \sum_{n=1}^{\infty}\int_{0} ^{t}f(s)e_{n}db^{H}_{n}(s)&=&\sum_{n=1}^{\infty}\sum_{m=1}^{\infty}e_{m}\int_{0}^{t}\langle f(s)e_{n}, e_{m}\rangle_{\mathcal{H}}\, db_{n}^{H}(s)\\&=&\sum_{n=1}^{\infty}\sum_{m=1}^{\infty}e_{m}\int_{0}^{t}\langle K_{H}^{*}(f(s)e_{n}), e_{m}\rangle_{\mathcal{H}}\, db_{n}(s)\\&=&\sum_{n=1}^{\infty}\int_{0}^{t}K_{H}^{*}(f(s)e_{n})\, db_{n}(s).
        \end{eqnarray*}

\subsection{Smoothing effect}
Let ${\rm e}^{t\Delta}$ be the linear heat semi-group defined by ${\rm e}^{t\Delta}\,\varphi=G_t\star\varphi, t>0,$ where $G_t$ is the heat kernel given by
$$
G_t(x)=\left(4\pi t\right)^{-N/2}\,{\rm e}^{-\frac{|x|^2}{4t}},\;\;\;t>0,\;\;\;x\in\R^N.
$$
Let, for $\gamma\geq 0$, ${\mathbf S}_{\gamma}$ be defined as
\begin{equation}
\label{S}
{\mathbf S}_{\gamma}(t)\varphi={\rm e}^{t\Delta}\left(|\cdot|^{-\gamma}\varphi\right),\;\;\;t>0.
\end{equation}
To treat the nonlinear term in \eqref{main}, we use the following key estimate proved in \cite{BTW}.
\begin{proposition}
\label{Key}
Let $N\geq 1$ and $0< \gamma<N$. For $1<q_1, q_2\leq \infty$ such that
\begin{equation}
\label{q12}
\frac{1}{q_2}<\frac{\gamma}{N}+\frac{1}{q_1}<1,
\end{equation}
we have
\begin{equation}
\label{SS}
\|{\mathbf S}_{\gamma}(t)\varphi\|_{q_2}\leq {\mathtt C}_0 \,t^{-\frac{N}{2}\left(\frac{1}{q_1}-\frac{1}{q_2}\right)-\frac{\gamma}{2}}\,\|\varphi\|_{q_1},
\end{equation}
where ${\mathtt C}_0 $ is a constant depending only on $N, \gamma, q_1$ and $q_2$.
\end{proposition}
\begin{remark}
For $\gamma=0$, the estimate \eqref{SS} holds under the assumption $1\leq q_1\leq q_2\leq\infty$. This is unlike to the case $\gamma>0$, where \eqref{q12} enable us to take $q_2<q_1$.
\end{remark}
\section{Proof of the main result}
This section is devoted to the proof of Theorem \ref{LWP}. First, we consider the linear Cauchy problem
\begin{equation}
\label{Lin}
\left\{
\begin{matrix}
\partial_t{\mathbf Z}(t)-\Delta {\mathbf Z}(t)=\partial_t B^H(t),\quad t>0,\\
{\mathbf Z}(0)= 0.\\
\end{matrix}
\right.
\end{equation}
The mild solution of \eqref{Lin} is given by
$$
{\mathbf Z}(t)=\int_0^t\,{\rm e}^{(t-s)\Delta}\,d\,B^H (s).
$$
Thanks to \eqref{H} and according to \cite{CMO}, we know that the mild solution ${\mathbf Z}$ belongs to $C([0,T]; L^q)$ for any $T>0$. For $\sigma>0$ and $1<r<\infty$ to be fixed later, define
\begin{equation}
\label{K}
{\mathbf K}(T)=\sup_{0\leq t\leq T}\|{\mathbf Z}(t)\|_q+\sup_{0\leq t\leq T}\,\left(t^{\sigma}\,\|{\mathbf Z}(t)\|_r\right).
\end{equation}
Now we are ready to give the detailed proof of Theorem \ref{LWP}.
\begin{proof}[{Proof of Theorem \ref{LWP}}]
We solve \eqref{main} via the integral form \eqref{integral}. Clearly, the mild solution $u$ can be seen as a fixed point of
\begin{equation}
\label{phi}
\Phi(u)(t)= {\rm e}^{t\Delta}u_{0}+\int_{0}^{t}\,{\mathbf S}_{\gamma}(t-s)\left(|u(s)|^{p-1}u(s)\right)\,ds+{\mathbf Z}(t).
\end{equation}
Let $r, \sigma$ be given by
\begin{equation}
\label{r}
\frac{1}{r}=\frac{1}{2qp}+\frac{1}{2q}-\frac{\gamma}{2N p},
\end{equation}
and
\begin{equation}
\label{sigma}
\sigma=\frac{N}{2}\left(\frac{1}{q}-\frac{1}{r}\right).
\end{equation}

Using \eqref{gamma} and \eqref{Q}, we see that $r>p>q>1$ and
\begin{equation}
\label{pqr}
\frac{1}{r}-\frac{1}{qp}+\frac{\gamma}{N p}<\frac{2}{Np},\quad \frac{1}{q}<\frac{\gamma}{N}+\frac{p}{r}<1,\quad 0<\sigma<\frac{1}{p}.
\end{equation}
For $M, T>0$ define
$$
{\mathbf X}=\Big\{\, u\in C([0,T]; L^q)\,\cap\, C((0,T]; L^r);\;\;\; \sup_{0\leq t\leq T}\,\|u(t)\|_q\leq M \;\;\;\mbox{and}\;\;\; \sup_{0< t\leq T}\,t^{\sigma}\|u(t)\|_r\leq M\,\Big\},
$$
endowed with the metric
$$
d(u,v)=\sup_{0\leq t\leq T}\,\|u(t)-v(t)\|_q+\sup_{0< t\leq T}\,t^{\sigma}\|u(t)-v(t)\|_r.
$$
Clearly $\left({\mathbf X},d\right)$ is a complete metric space. The rest of the proof uses a contraction mapping argument in the complete metric space $\left({\mathbf X},d\right)$.\\
First let us show that $\Phi({\mathbf X})\subset {\mathbf X}$ for suitable $M$ and $T$. Given $u\in {\mathbf X}$, we have
\begin{eqnarray*}
\|\Phi(u)(t)\|_q&\leq& \|{\rm e}^{t\Delta}u_{0}\|_q+\int_{0}^{t}\,\|{\mathbf S}_{\gamma}(t-s)\left(|u(s)|^{p-1}u(s)\right)\|_q\,ds+\|{\mathbf Z}(t)\|_q\\
&\leq& \|u_{0}\|_q+\|{\mathbf Z}(t)\|_q+{\mathtt C}_0\,\int_0^t\,(t-s)^{-\frac{N}{2}(\frac{p}{r}-\frac{1}{q})-\frac{\gamma}{2}}\,\||u(s)|^{p-1}u(s)\|_{\frac{r}{p}}\,ds\\
&\leq&\|u_{0}\|_q+\|{\mathbf Z}(t)\|_q+{\mathtt C}_0\,M^{p}\,\int_0^t\,(t-s)^{-\frac{N}{2}(\frac{p}{r}-\frac{1}{q})-\frac{\gamma}{2}}\,s^{-p\sigma}\,ds,
\end{eqnarray*}
where we have used \eqref{SS} with $q_1=\frac{r}{p}$, $q_2=q$. Note that, for $0\leq t\leq T$, we have $\|{\mathbf Z}(t)\|_q\leq {\mathbf K}(T)$ where ${\mathbf K}$ is given by \eqref{K}. In addition, remark that
$$
\int_0^t\,(t-s)^{-\frac{N}{2}(\frac{p}{r}-\frac{1}{q})-\frac{\gamma}{2}}\,s^{-p\sigma}\,ds=t^{1-p\sigma-\frac{N}{2}(\frac{p}{r}-\frac{1}{q})-\frac{\gamma}{2}}\,
{\mathcal B}\left(1-\frac{N}{2}(\frac{p}{r}-\frac{1}{q})-\frac{\gamma}{2}, 1-p\sigma\right),
$$
where the beta function ${\mathcal B}$ is given by
$$
{\mathcal{B}}(x,y)=\int_0^1\,\tau^{x-1}(1-\tau)^{y-1}d\tau,\quad x,\; y>0.
$$
Using \eqref{pqr}, we see that $1-p\sigma>0$ and $1-\frac{N}{2}(\frac{p}{r}-\frac{1}{q})-\frac{\gamma}{2}>0$. Moreover, using \eqref{sigma} and \eqref{Q}, it follows that
$$
1-p\sigma-\frac{N}{2}(\frac{p}{r}-\frac{1}{q})-\frac{\gamma}{2}=\frac{2-\gamma}{2q}\left(q-\frac{N(p-1)}{2-\gamma}\right)>0.
$$
Set $\alpha=\frac{2-\gamma}{2q}\left(q-\frac{N(p-1)}{2-\gamma}\right)$ and ${\mathtt C}_1={\mathcal B}\left(1-\frac{N}{2}(\frac{p}{r}-\frac{1}{q})-\frac{\gamma}{2}, 1-p\sigma\right)$, we end up with
\begin{equation}
\label{stab}
\|\Phi(u)(t)\|_q\leq \|u_0\|_q+{\mathbf K}(T)+ {\mathtt C}_0\,{\mathtt C}_1\,M^{p}\,T^{\alpha}.
\end{equation}
Similarly, using the smoothing effect of the heat semi-group and Proposition \ref{Key}, we get
\begin{eqnarray*}
\|\Phi(u)(t)\|_r&\leq& \|{\rm e}^{t\Delta}u_{0}\|_r+\int_{0}^{t}\,\|{\mathbf S}_{\gamma}(t-s)\left(|u(s)|^{p-1}u(s)\right)\|_r\,ds+\|{\mathbf Z}(t)\|_r\\
&\leq& t^{-\sigma}\|u_{0}\|_q+\|{\mathbf Z}(t)\|_r+{\mathtt C}_0\,\int_0^t\,(t-s)^{-\frac{N(p-1)}{2r}-\frac{\gamma}{2}}\,\||u(s)|^{p-1}u(s)\|_{\frac{r}{p}}\,ds\\
&\leq&t^{-\sigma}\|u_{0}\|_q+\|{\mathbf Z}(t)\|_r+{\mathtt C}_0\,M^{p}\,\int_0^t\,(t-s)^{-\frac{N(p-1)}{2r}-\frac{\gamma}{2}}\,s^{-p\sigma}\,ds.
\end{eqnarray*}
It follows that
\begin{equation}
\label{stabb}
t^{\sigma}\,\|\Phi(u)(t)\|_r\leq \|u_0\|_q+{\mathbf K}(T)+{\mathtt C}_0\,{\mathtt C}_2\,M^{p}\,T^{\alpha},
\end{equation}
where ${\mathtt C}_2={\mathcal B}(1-\frac{N(p-1)}{2r}-\frac{\gamma}{2}, 1-p\sigma)$.
Choose $M>\|u_0\|_q$. We infer, by choosing $T$ small enough, that
\begin{equation}
\label{T1}
{\mathbf K}(T)+ {\mathtt C}_0\,{\mathtt C}_1\,M^{p}\,T^{\alpha}\leq M-\|u_0\|_q,
\end{equation}
and
\begin{equation}
\label{T2}
{\mathbf K}(T)+ {\mathtt C}_0\,{\mathtt C}_2\,M^{p}\,T^{\alpha}\leq M-\|u_0\|_q.
\end{equation}
This proves that $\Phi(u)\in {\mathbf X}$.\\
Let us prove now that $\Phi$ is a contraction. For $u,v \in {\mathbf X}$, we have by using similar arguments as above
\begin{equation*}
d(\Phi(u),\Phi(v))\leq\,{\mathtt C}_0\,({\mathtt C}_1+{\mathtt C}_2)\,M^{p-1}\,T^{\alpha}\,d(u,v)\leq \,\frac{1}{2}\,d(u,v),
\end{equation*}
provided that
\begin{equation}
\label{cont1}
{\mathtt C}_0\,({\mathtt C}_1+{\mathtt C}_2)\,M^{p-1}\,T^{\alpha}\,\leq\,\frac{1}{2}.
\end{equation}
Finally, for $M>\|u_0\|_q$ and $T$ satisfying \eqref{T1}-\eqref{T2}-\eqref{cont1}, we conclude that $\Phi$ is a contraction from ${\mathbf X}$ into itself. This finishes the proof of Theorem \ref{LWP}.
\end{proof}

\end{document}